\documentclass{amsart}

\usepackage{amssymb}

\theoremstyle{plain}

\numberwithin{equation}{section}

\newtheorem{theorem}[equation]{Theorem}

\newtheorem{hypotheses}[equation]{Hypotheses}

\newtheorem{lemma}[equation]{Lemma}

\newtheorem{example}[equation]{Example}

\newcommand{\Irr}{\operatorname{Irr}}

\newcommand{\Ker}{\operatorname{Ker}}

\newcommand{\Aut}{\operatorname{Aut}}

\newcommand{\dl}{\operatorname{dl}}

\newtheorem*{teoremaa}{Theorem A}

\newtheorem*{teoremab}{Theorem B}

\newtheorem*{teoremac}{Theorem C}
\newcommand{\core}{\operatorname{core}}

\theoremstyle{definition}

\newtheorem{definition}[equation]{Definition}

\begin{document}

	\title{Products of characters and derived length}

\author{Edith Adan-Bante}

\address{Department of Mathematics, 
University of Illinois at Urbana-Champaign,
Urbana, 61801}

\email{adanbant@math.uiuc.edu}


\thanks{This research was partially supported by grant NSF DMS 9970030}

\subjclass{20c15}

\date{2002}

\begin{abstract}
 Let $G$ be a finite solvable group and $\chi\in \Irr(G)$ be a faithful
character. We 
show that the derived length of $G$ is bounded by a linear 
function of the number of distinct irreducible 
constituents of  $\chi\overline{\chi}$. We also discuss 
other properties of the decomposition of $\chi\overline{\chi}$
into its irreducible constituents.
\end{abstract}

\maketitle\begin{section}{Introduction}
Let $G$ be a finite group. 
Denote by $\Irr(G)$ the set of 
irreducible complex characters of $G$. Let 
$1_G$ be the principal character of $G$. Denote by 
$[\Theta, \Phi]$ the inner product of the characters
$\Theta$ and $\Phi$ of $G$. Through this work, we will use
the notation of \cite{isaacs}.

Let $\chi \in \Irr(G)$. Define $\overline{\chi}(g)$ to be
the complex conjugate $\overline{\chi(g)}$ of $\chi(g)$ 
for all $g \in G$. Then $\overline{\chi}$ is also 
an irreducible complex character of $G$. Since the product of 
characters is a character, $\chi\overline{\chi}$ is a 
character of $G$. So it can be expressed as an integral linear
combination of irreducible characters. 
Now observe that 
	$$[\chi\overline{\chi}, 1_G]=[\chi, \chi]=1,$$
\noindent where the last equality holds 
 since $\chi \in \Irr(G)$.
Assume now that $\chi(1)>1$. Then the decomposition of
the character $\chi \overline{\chi}$ into its 
distinct irreducible 
constituents $1_G$, $ \alpha_1$, $\alpha_2, \ldots, \alpha_n$
has the form
\begin{equation}\label{suma2}
 \chi\overline{\chi}= 1_G + \sum_{i=1}^n a_i \alpha_i,
\end{equation}
\noindent where $n>0$ and 
$a_i>0$ is the multiplicity of 
$\alpha_i$.  

 Set $\eta(\chi)=n$, so that $\eta(\chi)$ 
is the number of distinct non-principal irreducible constituents
of $\chi\overline{\chi}$.
The number $\eta(\chi)$ carries information about the structure
of the group. For example, if $\eta(\chi)$ is an odd number,
 then the order of the group has to be an even number. 
To see this, notice that
$\chi\overline{\chi}$ is a real character. When $\eta(\chi)$ is odd, 
at least one of the irreducible characters $\alpha_i$ has to be
real. Then $G$ has a non-principal irreducible real character.
So the order of $G$ has to be even.
	
	The purpose of this work is to give some answers to the following 
questions:
	
	{\bf Question 1.}
Assume that we know $\eta(\chi)$ for some
$\chi \in \Irr(G)$.  What can we say about the structure of the 
group $G$ and  about the character $\chi$?

	{ \bf Question 2.} 
	Knowing the set $\{a_i \mid i=1, \dots, \eta(\chi)\}$, what can we say 
about the group $G$?

	Denote by $\dl(G)$ the derived length of the group $G$. 
	The main results of this work regarding the first question 
are the following:

\begin{teoremaa}
There exist constants $C$ and $D$ such that for any finite
solvable group $G$ and any irreducible character $\chi$
\begin{equation*}
\dl(G/\Ker(\chi)) \leq C \eta(\chi) +D.
\end{equation*}

\end{teoremaa}

\begin{teoremab}
Let $G$ be a finite solvable group and $\chi \in \Irr(G)$.
Then $\chi(1)$ has at most 
$\eta(\chi)$ distinct
prime divisors. 

If, in addition, $G$ is supersolvable and $\chi(1)>1$,
then  $\chi(1)$ is a product of at most $\eta(\chi)- 1$ primes.
\end{teoremab}

	The main result of this work regarding the second question
is

\begin{teoremac}
Assume that  $G$ is a finite solvable group and $\chi \in \Irr(G)$ 
with $\chi(1)>1$.
Let $\{\alpha_i \in \Irr(G)^
{\#}\mid i=1, \ldots , n \}$
be the set of 
non-principal irreducible constituents of $\chi \overline{\chi}$.
If   
 $\, \Ker(\alpha_j)$ is maximal under inclusion among the subgroups 
$\, \Ker(\alpha_i)$, for $i=1, \ldots, n$, of $G$,
 then $ [ \chi \overline{\chi} , \alpha_j] =1$.
Thus $1\in \{ [\chi \overline{\chi} , \alpha_i] \mid 
 i=1, \ldots, n\}$.
\end{teoremac}

{\bf Notation.} Set  $V^{\#}=V \setminus \{0\}$
and $\Irr(G)^{\#}= \Irr(G)\setminus \{ 1_G\}$.
\end{section}

\begin{section}{Preliminaries}

\begin{definition}
Let $V$ be a finite  ${\bf F} G$-module for some finite field ${\bf F}$.
Then $ m(G,V)$ is the number of distinct sizes of  
orbits of $G$ on $V^{\#}$. 
\end{definition}

\begin{lemma} [Keller] \label{orbits}

There exist universal constants $C_1$ and $C_2$ such that for any
 finite solvable group $G$ acting faithfully and irreducibly on a 
finite vector space $V$ we have  

\begin{equation*}
\dl(G) \leq C_1 \log ( m(G, V)) +C_2.
\end{equation*}
\end{lemma}
\begin{proof}
See \cite{keller}.
\end{proof}     

\begin{definition}\label{horbits} We define the function 
\begin{equation*}
 h(n)= C_1 \log ( n ) +C_2
\end{equation*}
\noindent where $C_1$ and $C_2$ are  as in Lemma  \ref{orbits}.
\end{definition}
\end{section}

\begin{section}{The function $\eta(\chi)$}

Given a finite group $G$ and a character $\chi \in \Irr(G)$, we 
 define $\eta(\chi)$ as the number of non-principal irreducible
constituents of the product $\chi \overline{\chi}$.  We  give
examples showing 
 that there is no relation between induction of characters and
$\eta$.




\begin{example} 
If $\chi=\theta^G$ is induced from some $\theta \in \Irr(H)$, where
$H\leq G$,  then we need not have     
$\eta(\chi) \geq \eta(\theta)$. 
\end{example}
\begin{proof}
Let $E$ be an extra-special group of exponent $p$ and order $p^3$,
for some odd prime $p$. 
Let $a \in \Aut(E)$ be an element of prime order $q$ that divides
$p-1$. Assume that
$a$ acts fixed point free on $E$.

	Set $G= <a>E$. Let $\theta\in \Irr(E)$ be a non-linear character.
Since $a$ acts fixed point free, we have that $\theta^G=\chi \in \Irr(G)$.

	Observe that $G$ has $q$ linear characters, 
namely the irreducible characters
of $G/E$. Also $G$ has $\frac{p^2-1}{q}$ irreducible characters of degree
$q$, the characters that are induced from linear non-principal
 characters of $E$.
And finally there are $\frac{p-1}{q}$ irreducible characters
 of degree $pq$.
We conclude that $G$ has $q + \frac{p^2-1}{q}+ \frac{p-1}{q}$
distinct irreducible characters.
Thus $\eta({\chi}) \leq q-1 + \frac{p^2-1}{q}+ \frac{p-1}{q}.$

	We can check that
 $$q-1 + \frac{p^2-1}{q}+ \frac{p-1}{q}<p^2-1.$$
\noindent 
Observe that $\theta\overline{\theta}=(1_{{\bf Z}(E)})^E$.
Thus $\eta(\theta)=p^2-1 > \eta(\chi)$.
\end{proof}
\begin{example}
If $\chi=\theta^G$ is induced from some $\theta \in \Irr(H)$, 
where
$H\leq G$,  then we need not have     
$\eta(\chi) \leq \eta(\theta)$. 
\end{example}
\begin{proof}
Let $G$ be an extra-special group. Let $\chi\in \Irr(G)$ be a non-linear
character. Let $\theta $ be a linear character of some 
subgroup $H$ of $G$  such that $\chi = \theta^G$.
Then $\eta(\chi)> \eta(\theta)=0$.
\end{proof}

\end{section}

\begin{section}{Proof of Theorem C}

	Let $G$ be a finite group and $\chi \in \Irr(G)$.
 Consider the expression \eqref{suma2} for $\chi\overline{\chi}$. 
We will see in this section that if $G$ is 
solvable, then $1 \in \{a_i\}$. That may not be true in general.
For example, consider $A_6$, the alternating group on 
6 letters,  and $\chi_5 \in \Irr(A_6)$ with 
$\chi_5(1)=10$. Using the notation
of page 289 of \cite{isaacs},   we can check that
$$\chi_5\overline{\chi_5} = \chi_1 + 2\chi_2 + 
2\chi_3+ 3\chi_4 +2\chi_5 + 2\chi_6 + 2\chi_7.$$
\noindent Thus $\{a_i\}=\{2,3\}$. 

\begin{lemma} \label{basico1}
Let $L$ and $N$ be normal subgroups of $G$ such that
$L/N$ is an abelian chief factor of $G$. Let $\theta \in \Irr(L)$
be a $G$-invariant character. Then the restriction 
$\theta_N$ is reducible if and 
only if 
\begin{equation}\label{0outside}
	\theta(g)=0\text{ for all } g \in L \setminus N.
\end{equation}
Also if $\theta_N$ is reducible, then 
 \begin{equation}\label{principalcharacter}
\theta\overline{\theta}= (1_N)^L + \Phi
\end{equation}
\noindent where 
$\Phi$ is either the zero function or a 
character of $L$, and $[\Phi_N, 1_N]=0$.
\end{lemma}
\begin{proof}
Let $\varphi\in \Irr(N)$ be a character such that 
$[\varphi, \theta_N]\neq 0$.
If  $\theta_N$ is reducible, by Theorem 6.18 of \cite{isaacs}
we have that either $\theta_N= e \varphi$, where $e^2=|L:N|$, or
$\theta= \varphi^L$. If 
$\theta_N= e \varphi$, where $e^2=|L:N|$, 
by Exercise 6.3 of \cite{isaacs}
we have that 
$\theta$ vanishes on $L \setminus N$. 
 If $\theta=\varphi^L$, since $N$ is a normal subgroup of 
$L$ we have that $\theta(g)=0$ for all $g \in L \setminus N$.
Thus \eqref{0outside} holds. 

Now assume that  \eqref{0outside} holds. Then

\begin{equation*} 
\begin{split}
[\theta_N, \theta_N] & =\frac{1}{|N| } \sum_{g \in N }\theta(g)
\overline{\theta(g)}  \\ 
& = \frac{1}{|N| } \sum_{g \in L }\theta(g)
\overline{\theta(g)} \qquad  \text{ by \eqref{0outside}}  \\ 
&= \frac{1}{|N|} |L| [\theta, \theta]= \frac{|L|}{|N|},
\end{split}
\end{equation*}
\noindent where the last equality holds since $\theta\in \Irr(L)$. 
Because $|L|/|N|>1$, it follows that $\theta_N$ is a reducible
character. 

For any $\gamma \in \Irr(L/N)$ we have that
{\allowdisplaybreaks 
\begin{equation*}
\begin{split}
[\theta\overline{\theta}, \gamma] & = [\theta, \theta\gamma]\\
&=\frac{1}{|L| } \sum_{g \in L }\theta(g)
\overline{\theta(g)} \overline{\gamma(g)}\\
&=\frac{1}{|L| } [\sum_{g \in L \setminus N} \theta(g)\overline{\theta(g)}
\overline{\gamma(g)} +
\sum_{g\in N} \theta(g)\overline{\theta(g)}\overline{\gamma(g)}]\\
&= \frac{1}{|L| } [\sum_{g\in L \setminus N} \theta(g)\overline{\theta(g)} +
\sum_{g\in N} \theta(g)\overline{\theta(g)\gamma(g)}] 
 \qquad \  \text{ by \eqref{0outside}} \\
&= \frac{1}{|L| } [\sum_{g\in L \setminus N} \theta(g)\overline{\theta(g)} +
\sum_{g\in N} \theta(g)\overline{\theta(g)}] \ \  \text{ 
since $\Ker(\gamma)\geq N$ and $\gamma(1)=1$} \\
&= [\theta,\theta]=1.
\end{split}
\end{equation*}
}
Thus  \eqref{principalcharacter} follows.
\end{proof}
\begin{lemma}\label{basicon}
 Let $G$ be a finite solvable group and $\chi\in \Irr(G)$.
Let $\{\alpha_i\mid i=1, \ldots, \eta(\chi)\}$ be the set
of non-principal 
irreducible constituents of the product
$\chi \overline{\chi}$.
 Let  $N$ be
a normal 
subgroup of $G$. Then $\chi_N\in \Irr(N)$ if and only if 
$N \not\leq\Ker(\alpha_i)$ for $i=1, \ldots,
\eta(\chi)$.
\end{lemma} 
\begin{proof} Observe that 
{\allowdisplaybreaks
\begin{equation*}
\begin{split}
[\chi_N, \chi_N] & = [\chi_N \overline{\chi}_N, 1_N]\\
	& = [(1_G+\sum_{i=1}^n a_i  \alpha_i)_N, 1_N] \ 
 \qquad \mbox{ by \eqref{suma2}}\\
	&= [1_N+\sum_{i=1}^n a_i (\alpha_i)_N, 1_N]\\
 	& = [1_N, 1_N]+\sum_{i=1}^n a_i [ (\alpha_i)_N, 1_N]\\
	&= 1+ \sum_{i=1}^n a_i [ (\alpha_i)_N, 1_N].
\end{split}
\end{equation*}
}
Thus $[\chi_N, \chi_N]=1$ if and only if 
$\sum_{i=1}^n a_i [ (\alpha_i)_N, 1_N]=0$. Since $a_i>0$ for $i=1, \ldots
,n$, we have 
$[\chi_N, \chi_N]=1$ if and only if 
$[ (\alpha_i)_N, 1_N]=0$ for  all $i$. 
Since $[ (\alpha_i)_N, 1_N]=0$ if and only if $N \not\leq \Ker(\alpha_i)$,
the result follows.
\end{proof}


\begin{proof}[Proof of Theorem C]
Set $N=\Ker(\alpha_j)$.
Let $L$ be a normal subgroup of $G$ such that $L/N$ is 
a chief factor of $G$. Since $N=\Ker(\alpha_j)\not\leq \Ker(\alpha_i)$
for $i=1, \ldots, n$, we have
 $L \not\leq \Ker(\alpha_i)$ for 
$i=1, \ldots ,n$.  By Lemma \ref{basicon}
we have that $\chi_L \in \Irr(L)$. Set $\theta= \chi_L$.
Since $N=\Ker(\alpha_j)$, we have that  
$[(\alpha_j)_N, 1_N]=\alpha_j(1)$. Thus 
$$[\chi_N, \chi_N]= [(\chi\overline{\chi})_N , 1_N]\geq 
1+a_j\alpha_j(1)>1.$$
\noindent Therefore $\chi_N$ is reducible.
 By Lemma \ref{basico1} we have that
\begin{equation}\label{equationprincipal} 
(\chi\overline{\chi})_L = \theta\overline{\theta}= 1_{N}^L +\Phi,
\end{equation}
\noindent where $\Phi$ is either the zero function or a 
character of $L$ and 
$[\Phi_N, 1_N]=0$. Also, by \eqref{suma2}
we have that
	\begin{equation*} 
(\chi\overline{\chi})_L = 1_L + \sum_{i=1}^n a_i(\alpha_i)_L.
\end{equation*}
Let $\gamma \in \Irr(L / {\Ker(\alpha_j)})$ be such that 
$[(\alpha_j)_L, \gamma]\neq 0$.
Then 
{\allowdisplaybreaks
\begin{equation*}
0< a_j [(\alpha_j)_L, \gamma] 
= [(a_j \alpha_j)_L, \gamma]\leq [ (\chi\overline{\chi})_L, \gamma]=1,
\end{equation*}
}
\noindent 
where the last equality follows  from \eqref{equationprincipal}.
Therefore  $a_j=1$.
 
Since there is  some $j\in \{1, \ldots ,n\}$ such that 
$\Ker(\alpha_j)$ is maximal among the  
$\Ker(\alpha_i)$ for all $i$, the last part of Theorem C
follows from that.
\end{proof}
\end{section}

\begin{section}{Proof of Theorem B}

\begin{lemma}\label{center}
Assume $G$ is a finite group and $\chi \in \Irr(G)$ is a faithful character. 
Let
$\{\alpha_i \in \Irr(G)^{\#}\mid i=1, \ldots , n \}$
be the set of 
non-principal irreducible constituents of $\chi \overline{\chi}$.
 Then
\begin{equation*}
{\bf Z}(G)= \bigcap_{i=1}^n \Ker(\alpha_i).
\end{equation*}
\end{lemma}
\begin{proof}
By Lemma 2.21 of \cite{isaacs}, 
$$\Ker(\chi\overline{\chi})= \Ker(1_G) \bigcap_{i=1}^n \Ker(\alpha_i)=
\bigcap_{i=1}^n \Ker(\alpha_i).$$
Since $(\chi\overline{\chi})(g)= \chi^2(1)$ if and only if
$g\in {\bf Z}(\chi)$, it follows that $\Ker(\chi\overline{\chi})=
{\bf Z}(G)$, and the result follows.
\end{proof}

 \begin{definition}\label{maximalchain}
 Let $G$ be a group and $L$ be a subgroup of $G$. We 
say that 
	$$(N, \theta)\leq (L, \phi)$$
\noindent if $N\leq L$, 
$\phi \in \Irr(L)$, $\theta \in \Irr(N)$
 and $[\phi_N,\theta]\neq 0$. We say that 
$$(N, \theta)< (L, \phi)$$
\noindent if $N < L$, 
$\phi \in \Irr(L)$, $\theta \in \Irr(N)$
 and $[\phi_N,\theta]\neq 0$.

Let $X$ be a family of normal subgroups of $G$ with
$G \in X$.  
We say that a chain
\begin{equation*} 
	(N_0, \theta_0) >
(N_1, \theta_1)> (N_2, \theta_2)>  \cdots > (N_k , \theta_k),
\end{equation*}
\noindent where $N_0=G$ and $\chi=\theta_0$, 
 is an {\bf $(X, \chi)$-reducing chain} if $N_i \in X$ and 
$(\theta_i)_{N_{i+1}}$ is reducible for $i=0,\ldots, k$. 

	We say that the above
 chain is a {\bf maximal  $(X, \chi)$-reducing chain}
if it is a  $(X, \chi)$-reducing chain with the following two
properties:
	
	(i) For any $i$ with $0<i\leq k$, 
the group 
$N_i$ is a maximal subgroup in the set
$$\{ M \in X \mid M 
\leq N_{i-1}\mbox{ and } (\theta_{i-1})_{M}\mbox{ is reducible}\}.$$  
	
	(ii) For any $M\in X$ such that 
$M< N_k$, the restriction $(\theta_k)_M$ is
irreducible.
 \end{definition}

{\bf Remark.} Given  a family $X$ of normal subgroups of $G$
with $G \in X$
and given $\chi\in \Irr(G)$, 
 there is always an $(X,\chi)$-reducing chain, and
a maximal  $(X, \chi)$-reducing chain.
In fact   $(G,\chi)$ is already an $(X,\chi)$-reducing chain.
 We find a maximal reducing $(X,\chi)$ chain 
by induction. 
We start with $(N_0, \theta_0)=(G, \chi)$. If 
$(\theta_0)_M$ is irreducible for any
$M \in X$, then $(N_0, \theta_0)$ is our maximal  
$(X, \chi)$-reducing chain. Assume we
have found
$(N_{i-1}, \theta_{i-1})$ for some integer $i\geq 1$. 
If the set
$$\{ M \in X \mid M \leq N_{i-1}\mbox{ and }
 (\theta_{i-1})_{M}\mbox{ is reducible}\}$$
\noindent is non-empty, we choose $N_i$ to be any maximal element in 
this set, and $\theta_i$ to be any character in $\Irr(N_i)$
such that $[(\theta_{i-1})_{N_i}, \theta_i]>0$.
Otherwise we stop our chain with $k=i-1$. 

\begin{hypotheses}\label{hypo}
Assume $G$ is  a finite solvable group and $\chi \in \Irr(G)$ is a 
faithful character. Set $n=\eta(\chi)$.
Let
$\{\alpha_i\in \Irr(G)^{\#} \mid i=1, \ldots , n \}$
be the set of 
non-principal irreducible constituents of $\chi \overline{\chi}$.
Set
\begin{equation}
	\Omega= \{\, \bigcap_{i\in S} \Ker(\alpha_i) \mid
 S \subseteq \{1,2, \ldots,  \,n \} \},
\end{equation}
\noindent where $\cap_{i \in S} \Ker(\alpha_i)$ is taken to be $G$ when $S$
is empty.

Let 	\begin{equation*} 
	(G, \chi)= (N_0, \theta_0) >
(N_1, \theta_1)> \cdots >(N_k, \theta_k)
\end{equation*}
be a  maximal $(\Omega, \chi)$-reducing chain.	

\end{hypotheses}

\begin{lemma}\label{maximal}
 Assume Hypotheses \ref{hypo}. 
Then the  maximal $(\Omega, \chi)$-reducing chain	
has the following properties:

{\normalfont (a)} For any integer $i=1, 2, \ldots, k$ and any normal
subgroup $M$ of $G$  such that 
$N_{i}<  M \leq  N_{i-1}$ we have that 
\begin{equation}\label{m} 
(\theta_{i-1})_M \in \Irr(M). 
\end{equation} 
\indent {\normalfont (b)} $N_k$ is abelian.

{\normalfont	(c)} $k\leq n$. 

{\normalfont	(d)} If, in addition, $G$ is supersolvable, then 
$k\leq n-1$.
\end{lemma}
\begin{proof}{\bf  (a)}
If  $M \in \Omega$, then $(\theta_{i-1})_M$ has to be irreducible.
Otherwise $N_i$ is not a maximal element in $\Omega$ such that
$(\theta_{i-1})_{N_i}$ reduces, a contradiction with property
(i) in Definition \ref{maximalchain}.
 
So we may assume that $M$ is not an element of $\Omega$. Let $L$ be
minimal among all elements $K \in \Omega$ such that $M\leq K\leq N_{i-1}$. 
By property (i) in Definition \ref{maximalchain}
  we have that 
$$\phi=(\theta_{i-1})_L \in \Irr(L).$$
Observe	that
\begin{equation}\label{restrictionprocess}
1=[\phi,\phi] = [\phi \overline{\phi}, 1_L]
\leq [(\phi \overline{\phi})_M, 1_M]
= [\phi_M, \phi_M],
\end{equation}
\noindent where equality holds if and only if $\phi_M \in \Irr(M)$.

	Recall that  $[\chi_{N_{i-1}}, \theta_{i-1}] \neq 0$.
Thus $[\chi_L,\phi]\neq 0$.
Let $T$ be the stabilizer of 
$\phi$ in $G$ and $Y$ be a set of coset representatives of 
$T$ in $G$. Thus if $g,h \in Y$ and $g\neq h$,
 we have that 
 $\phi^g \neq \phi^h$ and therefore
$[\phi^g,\phi^h]=0$.
By Clifford
Theory we have that
$\chi_L= e \sum_{g\in Y} \phi^g$,
for some integer $e>0$. Thus
\begin{equation}\label{restrictionni}
[(\chi\overline{\chi})_L, 1_L] = [\chi_L,
\chi_L] 
 = [e\sum_{g\in Y} \phi^g, e\sum_{g\in Y} \phi^g]
= e^2 \sum_{g\in Y} [\phi^g,\phi^g] 
\end{equation}
	Since $\chi_M=(\chi_L)_M$, we have that
\begin{equation}\label{restrictionm}
\begin{split}
[(\chi\overline{\chi})_M, 1_M] &= [\chi_M,
\chi_M]\\
 &= e^2[\sum_{g\in Y} 
(\phi^g)_M, \sum_{g\in Y} (\phi^g)_M].
\end{split}
\end{equation}

	If $\phi_M\notin \Irr(M)$, 
then \eqref{restrictionprocess},
\eqref{restrictionni} and \eqref{restrictionm}  imply that
\begin{equation}\label{desigual}
[(\chi\overline{\chi})_L,1_L]<
 [(\chi\overline{\chi})_M,1_M].
\end{equation}
\noindent
By  \eqref{suma2} and  \eqref{desigual}
there  exists some $\alpha_j$ such that
 $\Ker(\alpha_j) \geq  M$ but
$\Ker(\alpha_j)\not\geq L$.
 Therefore $L\cap \Ker(\alpha_j)$ is a proper 
subset of $L$, contains $M$ and lies in $\Omega$.
This contradicts our
choice of $L$.
Thus  $(\theta_{i-1})_M=\phi_M \in \Irr(M)$.

{\bf (b)}
By Lemma \ref{center} we have that ${\bf Z}(G) \subseteq M$
for  any $M \in \Omega$. Thus $(\theta_k)_{{\bf Z}(G)}$ is irreducible
by property (ii) in Definition \ref{maximalchain}. 
That implies that 
	$\theta_k\in \Irr(N_k)$ is a linear character. Since 
$N_k$ is normal in $G$ 
and $[\chi_{N_k}, \theta_k]\neq 0$, all the irreducible
components of $\chi_{N_k}$ are linear. By hypothesis 
 $\chi\in \Irr(G)$ is a faithful character. Therefore 
$N_k$ must be abelian.

{\bf (c)}
This follows from the definition of $\Omega$
and the fact
that the set
$\{\Ker(\alpha_j)\}$ has at most $n$ elements.

{\bf (d)}
Suppose that $N_k={\bf Z}(G)$. Let $L/N_k$ be a chief factor of $G$ with 
$L\leq N_{k-1}$. Since $G$ is supersolvable, $L/N_k$ is cyclic 
of prime order. Observe that $L$ is abelian because
it has a central subgroup $N_k$ with a cyclic factor group
$L/N_k$. So  $\theta_k$ extends to $L$. 
By (a) we have that $(\theta_{k-1})_L \in \Irr(L)$. Thus 
$(\theta_{k-1})_{N_k}=\theta_k$. That can not be by Definition 
\ref{maximalchain} (i). We conclude that  $N_k \neq {\bf Z}(G)$.

 Since  $N_k \neq {\bf Z}(G)=\bigcap_{i=1}^n \Ker (\alpha_i)$
and $\{\Ker(\alpha_i) \mid i=1, 2,\ldots n\}$ has at most
$n$ elements,
 we must
have that $k\leq n-1$.
\end{proof}

Theorem B  is an application of Lemma \ref{maximal}.

\begin{proof}[Proof of Theorem B]
Working with the group $G/ \Ker(\chi)$, 
by induction on the order of $G$
we can assume that $\Ker(\chi)=1$. Let
\begin{equation*} 
	(G, \chi)=(N_0, \theta_0)>
(N_1, \theta_1)> \cdots > (N_k, \theta_k)
\end{equation*}
be a maximal $(\Omega, \chi)$-reducing chain.
For each $i=1, 2, \ldots, k$, let  $L_i$ be a
 normal subgroup of $G$ such that 
$L_i/ N_i$ is a chief factor of $G$ and $L_i \leq N_{i-1}$. 

By Lemma \ref{maximal} we have that
 $(\theta_{i-1})_{L_{i}} \in \Irr(L_{i})$.
Since $L_i/N_i$ is an elementary abelian $p_i$-group for some prime $p_i$,
 we have 
$$\theta_{i-1}(1)=\theta_i(1)  {p_i}^{m_i}$$
\noindent for some integer  $m_i\geq 1$. Here $m_i=1$ 
in the case that $G$ is supersolvable. By Lemma \ref{maximal} (b), we have
that $\theta_k(1)=1$. 
By Lemma \ref{maximal} (c),  $k\leq n$. We conclude that
$\chi(1)$ has at most $k\leq n$ distinct prime divisors.

If $G$ is supersolvable, by Lemma \ref{maximal} (d) we have $k \leq n-1$.
 Thus $\chi(1)$ has at
most $n-1$ prime divisors.
\end{proof}

\end{section}

\begin{section}{Proof of Theorem A}

\begin{hypotheses}\label{hypo2}
Assume Hypotheses \ref{hypo}. For each $i$, let $L_i/N_i$ be a chief factor of
$G$ where $L_i\leq N_{i-1}$.
\end{hypotheses}

\begin{lemma}\label{uinthemiddle}
Assume Hypotheses \ref{hypo2}.
 There exists a subgroup $U$ of $L_{i}$ and a character
$\phi \in \Irr(U)$, such that 
\begin{equation}\label{inthemiddle}
(N_{i}, \theta_{i}) \leq ( U, \phi) <( L_{i}, \psi)
\end{equation}
\end{lemma}
\begin{proof}
Suppose that the lemma  is false.
Then for any $U$ and $\phi \in \Irr(U)$ such that
\eqref{inthemiddle} holds, we have that $(L_{i})_{\phi}=L_{i}$.
 Choose a chain
$$(N_{i}, \theta_{i})=
(U_s, \phi_s)<\cdots <(U_1,\phi_1)<(U_0,\phi_0)=(L_{i},\psi)$$
\noindent such that $|U_{j-1}:U_j|$ is a prime number for
 all $j=1,2,\ldots,s$. We can do that
since $L_{i}/N_{i}$ is an elementary abelian group. Since 
 $(L_{i})_{\phi_j} =L_{i}$  for all $j=1,2,\ldots ,s$, we have 
 $ (U_{j-1})_{\phi_j} =U_{j-1}$.
Since $|U_{j-1}:U_j|$ is a prime number, it follows that
$(\phi_{j-1})_{U_j}=\phi_j$ for $j=1,\ldots, s$.
But then $(\theta_{i-1})_{N_{i}}\in \Irr(N_{i})$, a contradiction
with Definition \ref{maximalchain} (i). Therefore there exist
$U< L_{i}$ and a character $\phi \in \Irr(U)$ such that
\eqref{inthemiddle} holds and $ (L_{i})_{\phi}\neq L_i$.

 Since $N_i\leq (L_i)_{\phi} <L_i$, and 
$L_{i}/N_{i}$ is an elementary abelian subgroup, 
the subgroup $ (L_i)_{\phi}$
is normal in $L_i$.  By Clifford Theory $\psi$ is induced 
from some character $\psi_{\phi} \in \Irr( (L_i)_{\phi})$.
Since $(L_i)_{\phi}$
is normal in $L_i$, and $(\psi_{\phi})^{L_i}=\psi$, 
 we have $\psi(g)=0$ for any $g \in L_i\setminus (L_i)_{\phi}$.
\end{proof}

\begin{lemma}\label{plusone1}
 Assume Hypotheses \ref{hypo2}. 
Let $r_i = |\{ \alpha_j | N_i \leq \Ker(\alpha_j) \mbox{ and } N_{i-1} 
\not\leq \Ker(\alpha_j)\}|$.
Then  we have 
\begin{equation*}
 \dl(N_{i-1}/ {\bf C}_{N_{i-1}}(L_{i}/N_{i})) \leq h(r_i),
\end{equation*}
\noindent where $h$ is as in Definition \ref{horbits}.
\end{lemma}

\begin{proof}
By Lemma \ref{maximal} (a), we have that 
\begin{equation}
\psi=(\theta_{i-1})_{L_{i}} \in \Irr(L_{i}).
\end{equation}

Let $V(\psi)$ be the ``vanishing-off subgroup of $\psi$'' (see page
200 of  \cite{isaacs}), 
the smallest subgroup  $V(\psi)$ of $L_i$ such that
$\psi$ vanishes on $L_i \setminus V(\psi)$. Since 
$\psi=(\theta_{i-1})_{L_{i}}$ and $\theta_{i-1} \in \Irr(N_{i-1})$,
 the subgroup
$V(\psi)$ is $N_{i-1}$-invariant. Therefore $N_{i} V(\psi)$ is a
normal subgroup of $N_{i-1}$. Let $U$ and $\phi \in \Irr(U)$
be  as in 
Lemma \ref{uinthemiddle}. Observe that $V( \psi) \leq  (L_{i})_{\phi}$
since for 
 all $g \in L_{i}\setminus (L_{i})_{\phi}$ we have that 
$\psi(g)=0$.
 Also observe that $ N_{i} \leq (L_{i})_{\phi}$. 
Thus  $N_{i} V(\psi) \leq (L_{i})_{\phi}$. Therefore 
$N_{i} V(\psi)$ is a proper subgroup of $L_{i}$.  
 
Let $M$ be a subgroup such that $N_{i}V(\psi) \leq M< L_i$ and 
 $L_{i}/M$ is a chief 
factor of $N_{i-1}$. So we have the following relations:
$$N_i \leq N_i V(\psi)\leq  M < L_i \leq N_{i-1}.$$
 Since $L_{i}$ is a  normal subgroup 
of $G$, the quotient $L_i/M^g$ is also a chief factor of $N_{i-1}$, for any
$g \in G$.  Hence for any $g \in G$
\begin{equation}\label{twooptions}
  M M^g =M \
\mbox{  or } \ 
MM^g=L_{i}.
\end{equation}
  Lemma \ref{basico1} gives us that 
\begin{equation*}
\psi\overline{\psi}= 1_M^{L_{i}} + \Phi,
\end{equation*}
\noindent where $\Phi$ is either $0$ or a character of $L_{i}$. Since 
$[(\chi)_{L_{i}}, \psi]\neq 0$, this implies that 
\begin{equation*} 
(\chi\overline{\chi})_{L_{i}}= 1_M^{L_{i}} + \Theta,
\end{equation*}
\noindent where $\Theta$ is either $0$ or a character of $L_{i}$.
This and \eqref{suma2}  imply that
\begin{equation*}
1_{L_{i}} + \sum_{j=1}^n a_j (\alpha_j)_{L_{i}}= 1_M^{L_{i}} + \Theta.
\end{equation*}
Thus 
\begin{equation*}
\Irr(L_{i}/ M)^{\#}= \bigcup_{j=1}^n \{
\gamma \in \Irr( L_i/ M)^{\#} \mid
\ [(\alpha_j)_{L_i}, \gamma]\neq 0\}.
\end{equation*}
Let $X=\{\alpha_j | [(\alpha_j)_{L_i},\gamma]\neq 0 \ 
\mbox{ for some } \gamma \in  \Irr(L_{i}/ M)^{\#}\}$.
Observe that $X$ is a subset of the set
\begin{equation*}
\{ \alpha_j \  |  \ N_i \leq \Ker(\alpha_j)\  \mbox{ and } N_{i-1} 
\not\leq \Ker(\alpha_j)\}.
\end{equation*}
\noindent Thus 
\begin{equation}\label{riandx}
|X| \leq r_i.
\end{equation}
Let $\gamma, \delta \in \Irr(L_{i}/ M)^{\#}$.
Suppose that $\gamma$ and $\delta$ lie below the same $\alpha_j\in X$, 
i.e.
$[(\alpha_j)_{L_{i}}, \gamma]\neq 0$ and 
$[(\alpha_j)_{L_{i}}, \delta]\neq 0$,
for some $j=1, \ldots, n$.  
Since $L_{i}$ is a normal subgroup of $G$ and  $\alpha_j \in \Irr(G)$, 
by Clifford theory there exists $g \in G$ such that 
$\gamma^g=\delta$. By definition we have that $M\leq \Ker(\delta)$.
Observe that 
$$M^g\leq  (\Ker(\gamma))^g=\Ker(\gamma^g)$$
Since $\gamma^g=\delta$, we have $M M^g\leq \Ker (\delta)$.
By \eqref{twooptions}    
we have that $M^g=M$, i.e. $g \in N_G(M)$. We conclude that
 $\gamma$ and $\delta$ lie below the same $\alpha_j$ if and only if
$\gamma^g=\delta$ for some $g \in N_G(M)$, i.e the set 
$\{\gamma \in \Irr( L_i/ M)^{\#}\mid
 [(\alpha_j)_L, \gamma]\neq 0\}$ is an 
$N_G(M)$-orbit in $\Irr(L_i / M)^{\#}$.
Set $H=N_G(M)$. Each $H$-orbit in $\Irr(L_i/M)^{\#}$ lies under
at least one character $\alpha_j$ in $X$, and any 
each $\alpha_j$ lies over a single $H$-orbit $\Irr(L_i/M)$. 
Hence 
$H$ acts on $\Irr(L_{i}/M)^{\#}$ with at most $|X|$ orbits. 
By \eqref{riandx} we conclude that 
 $H$ acts on $\Irr(L_{i}/M)^{\#}$ with at most $r_i$ orbits.
By Lemma \ref{orbits} we have that
\begin{equation*}
\dl(H/{\bf C}_H(L_{i}/ M))\leq h(r_i).
\end{equation*}
	Since $N_{i-1} \leq H=N_G(M)$ and 
${\bf C}_H(L_{i}/ M)\cap N_{i-1}= {\bf C}_{N_{i-1}} (L_{i}/ M)$,
we have 
\begin{equation*}
\dl(N_{i-1}/{\bf C}_{N_{i-1}} (L_{i}/ M))\leq h(r_i).
\end{equation*}
	For any $g\in G$, we can check that 
\begin{equation*} 
	(G, \chi)= (N_0, (\theta_0)^g) >
(N_1, (\theta_1)^g)> \cdots > (N_k, (\theta_k)^g)
\end{equation*}
\noindent is a maximal $(G,\Omega)$-reducing chain. Thus, as 
before we can conclude that 
\begin{equation}\label{hwithrespect}
\dl(N_{i-1}/{\bf C}_{N_{i-1}} (L_{i}/ M^g))\leq h(r_i).
\end{equation}
Since $L_{i}/N_{i}$ is a chief factor of $G$ and 
$N_{i}\leq  M < L_{i}$,
we have that 
$$\core_G(M)=\bigcap_{g\in G} M^g= N_{i}.$$
Therefore
\begin{equation}\label{intersectingall}
\bigcap_{g\in G} { \bf C}_{N_{i-1}}(L_{i}/ M^g)= 
{ \bf C}_{N_{i-1}}(L_{i}/N_{i}).
\end{equation}
Observe that 
 the lemma follows from 
\eqref{hwithrespect} and  \eqref{intersectingall}.
\end{proof}
\begin{lemma}\label{plusone2}
 Assume Hypotheses \ref{hypo}. 
\begin{equation*}
\dl(N_{i-1}/N_{i})\leq \dl(N_{i-1}/ {\bf C}_{N_{i-1}}(L_{i}/N_{i})) + 1.
\end{equation*}
\end{lemma}
\begin{proof}
Set $C= {\bf C}_{N_{i-1}}(L_{i}/N_{i})$. Observe that 
$L_{i} \leq C$ and that $C$ is a normal subgroup of $G$.
 We want to prove that $C/N_i$ is abelian. We may assume
that $C>L_i$.
Observe that if $U$ is a group 
and  
$N_{i} \leq U \leq L_{i}$, 
then $U$ is normal in $C$. 
By Lemma \ref{uinthemiddle},
 there exist $U$ and $\phi \in \Irr(U)$, where
\begin{equation}\label{inthemiddlem}
(N_{i}, \theta_{i}) \leq ( U, \phi) <(L_i, \psi)
\end{equation}
\noindent and $(L_i)_{\phi} < L_i$. 
In particular we have that $C_{\phi}\neq C$. Since 
 $(\theta_{i-1})_{L_i} =\psi \in \Irr( L_i)$ and $U< L_i\leq C\leq N_{i-1}$,
 we have that 
$(\theta_{i-1})_C \in \Irr(C)$ and $(\theta_{i-1})_C$ lies above
$\phi$. By Clifford Theory, there exists $\zeta \in \Irr(C_{\phi})$ such
that $\zeta^C= (\theta_{i-1})_C$. Since $(\zeta^C)_{L_{i}} \in \Irr(L_i)$, 
we have that 
$C=C_{\phi} L_{i}$ (see Exercise 5.7 of \cite{isaacs}).
Observe that  $C_{\phi}$ is normal in $C$ since
$L_i/N_i$ is central in $C= {\bf C}_{N_{i-1}}(L_{i}/N_{i})$.
 Since $L_{i}/N_{i}$ is abelian, so is $C/C_{\phi}$.
Since $C$ is normal in $G$, for any $g \in G$ we have that 
$C/C_{\phi}^g$ is abelian.

Since $(\theta_{i-1})_{C} \in \Irr(C)$, while
 $[(\theta_{i-1})_U, \phi]\neq 0$
and $(L_i)_{\phi} < L_i$, we have
that $(\theta_{i-1})_{C_{\phi}}$ is a reducible character.
 Set $P=\bigcap_{g\in G} C_{\phi}^g$. Observe that $P$   
 is a normal subgroup of  $G$ with $N_i \leq P < N_{i-1}$.
Observe also that 
$(\theta_{i-1})_P$ is reducible since $P\leq C_{\phi}$.
 By Lemma \ref{maximal} (a), we
have that $P=N_{i}$.
Therefore  $C/ N_{i}$ is abelian and the lemma follows.
\end{proof}

\begin{lemma}\label{plusone}
 Assume Hypotheses \ref{hypo2}. Then
	\begin{equation*}
\dl(N_{i-1}/N_{i})\leq h(r_i) + 1.
\end{equation*}
\end{lemma}
\begin{proof}
It follows from Lemmas \ref{plusone1} and \ref{plusone2}
\end{proof}

\begin{lemma}\label{controlingfunction}
Let $n>1$ be an integer. Set  ${\bf N}=\{1,2,\ldots\}$.
Define
{\allowdisplaybreaks
\begin{equation}\label{dotsproduct}
p (n) = max \{ n_1\cdot n_2 \cdot \ldots \cdot n_s \ 
|\  n_1, n_2, \ldots, n_s \in {\bf N}
\mbox{ and } n_1 +n_2+\ldots+ n_s=n\}  
\end{equation}
}
Then $$p(n+1)\leq 2p(n).$$ Therefore
\begin{equation}\label{functionpn}
	p(n) \leq 2^{n-1}.
\end{equation}
\end{lemma}
\begin{proof} Observe that $n\leq p(n)$ since 
we can take $s=1$ and $n_1=n$ in \eqref{dotsproduct}. 
Thus if  $p(n+1)= m_1\cdot m_2 \cdot\ldots \cdot m_t$, where 
$m_1, m_2, \ldots, m_t$ are non-zero positive integers and 
$m_1+ m_2+  \ldots +m_t=n+1$, then  $m_i>1$ for some 
$i\in \{1, \ldots, t\}$.
 Assume that
$m_1\geq 2$. Then $m_1 -1\geq 1$, 
$ (m_1-1)+ m_2+  \ldots +m_t=n$. By definition we have that 
$	( m_1-1)\cdot m_2\cdot \ldots\cdot m_t \leq p(n).$
Thus 
{\allowdisplaybreaks
\begin{equation*}
\begin{split}
p(n+1)& =  m_1\cdot m_2 \cdot\ldots \cdot m_t\\
       & = ( m_1-1)\cdot m_2\cdot \ldots\cdot m_t+
 1\cdot m_2\cdot\ldots \cdot m_t\\
        & \leq p(n) +  1\cdot m_2\cdot\ldots \cdot m_t\\
		& \leq p(n) +  (m_1-1) \cdot m_2\cdot \ldots \cdot m_t\\
		&  \leq p(n)+p(n)=2p(n).
\end{split}
\end{equation*}
}
Since $p(2) = 2$, inequality \eqref{functionpn} follows. \end{proof}
\begin{proof}[Proof of  Theorem A]
Working with the group $G/ \Ker(\chi)$, 
by induction on the order of $G$
we can assume that $\Ker(\chi)=1$.
 So we may
assume Hypotheses \ref{hypo}.
Let
	\begin{equation*} 
	(G, \chi)= (N_0, \theta_0) >
(N_1, \theta_1)> \cdots > (N_k, \theta_k)
\end{equation*}
\noindent be a  maximal $(\Omega, \chi)$-reducing chain.
Set $n=\eta(\chi)$.
By Lemma \ref{maximal} (b) and (c), 
 we have that $N_k$ is abelian and $k\leq n$.
By Lemma \ref{plusone}, we have that, for $i=1,\ldots, k$,  
$$\dl(N_{i-1}/N_{i}) \leq h(r_i)+1, $$
\noindent where $r_i = |\{ \alpha_j | N_i \leq \Ker(\alpha_j)
\mbox{ and }  N_{i-1} 
\not\leq \Ker(\alpha_j)\}|$.
The definition of a maximal reducing 
chain and the definition of $r_i$ implies  that
\begin{equation}\label{countings}
r_1 + r_2 + \ldots + r_k\leq n. 
\end{equation}
By Lemma \ref{controlingfunction}
we have that 
\begin{equation*}
\prod_{i=1}^k r_i\leq 2^{n-1}.
\end{equation*}
Thus 
$$\dl(G)  \leq \sum_{i=1}^{k} \dl(N_{i-1}/N_i) + \dl(N_k)
	 \leq  \sum_{i=1}^{k} (h(r_i)+1) +1.$$
Since $h(r_i)= C_1 \log (r_i) + C_2$ by Definition 
\ref{horbits},
we have that 
\begin{equation*}
\begin{split} 
\dl(G)  \leq \  
\sum_{i=1}^{k} (C_1 \log (r_i) + C_2 +1)+1 = \
&  C_1 [\sum_{i=1}^k \log (r_i)] + (C_2+1)k +1\\
          \leq \  &  C_1 \  \log (\prod_{i=1}^k r_i) + (C_2+1)k +1.
	\end{split}
\end{equation*}
Let $s= \sum_{i=1}^k r_i$. By \eqref{countings} we have
$s\leq n$. By Lemma \ref{functionpn} we have that
	   $$ \prod_{i=1}^k r_i\leq 2^{s-1}\leq 2^{n-1}.$$        
Thus
 $$ \dl(G)\leq   C_1 \log (2^{n-1}) +  (C_2+1)k +1 
       	 \leq  (n-1) C_1 \log (2) +(C_2+1)n  + 1,$$ 
\noindent where the last inequality follows from $k\leq n$ (see Lemma
\ref{maximal} (c)). 
Set $C= C_1 \log (2) +C_2+1$ and $D= 1  + C_1 \log (2)$.
Then
	$$\dl(G) \leq C n + D.$$ 
\end{proof}

\begin{theorem}\label{supersolvable} 
Let $G$ be a supersolvable group. Let 
 $\chi \in \Irr(G)$ be such that $\chi(1)>1$. Then
	\begin{equation*}
	\dl(G/\Ker(\chi)) \leq 2\eta(\chi) -1.
\end{equation*}
\end{theorem}
\begin{proof}  
Working with the group $G/ \Ker(\chi)$, 
by induction on the order of $G$
we can assume that $\Ker(\chi)=1$.

 Let
	\begin{equation*} 
	(G, \chi)= (N_0, \theta_0) >
(N_1, \theta_1)> \cdots > (N_k, \theta_k)
\end{equation*}
\noindent be a  maximal $(\Omega, \chi)$-reducing chain.
Let $L_i/N_i$ be a chief factor of $G$, where $L_i \leq N_{i-1}$.
Since $G$ is a supersolvable group,  
$L_{i}/N_{i}$ is a cyclic group of prime order. Set 
$H=N_{i-1}/ {\bf C}_{N_{i-1}}(L_{i}/N_{i})$.
Observe that $H$  acts faithfully on $L_{i}/N_{i}$ as automorphisms. Since
$L_{i}/N_{i}$ is cyclic,  $H$ is abelian, i.e.

\begin{equation*}
 \dl(N_{i-1}/ {\bf C}_{N_{i-1}}(L_i/N_i)) \leq 1.
\end{equation*}

By Lemma \ref{plusone2} we conclude that
\begin{equation*}
\dl(N_{i-1} /N_i )\leq 2.
\end{equation*}

By Lemma \ref{maximal} (d) we have that $k\leq \eta(\chi)-1$. Also
$N_k$ is abelian by Lemma \ref{maximal} (b). 
Thus
$$\dl(G)\leq 2(\eta(\chi)-1)+1.$$
\end{proof}

\end{section}

{\bf Acknowledgment. } This is part of my Ph.D. Thesis.
I thank Professor Everett C. Dade, my adviser,
 and Professor
I. Martin Isaacs for their advise and suggestions. 
I would like to thank the mathematics department of the University 
of Wisconsin, at Madison, for 
their hospitality while I was visiting, and the mathematics department
of the University of  Illinois at Urbana-Champaign for their support.

\end{document}